\documentclass{amsart}
\usepackage{amssymb}
\usepackage{amsfonts}

\newtheorem{theorem}{Theorem}
\theoremstyle{plain}

\newtheorem{definition}{Definition}

\numberwithin{equation}{section}
\input{tcilatex}

\begin{document}
\title[]{Relaxed Elastic Line on an Oriented Surface in the Galilean Space}
\author{Tevfik \c{S}ahin}
\address{\textit{Amasya University,} \textit{Faculty of Sciences and Arts,} 
\textit{Department of Mathematics,} \textit{Amasya-Turkey}}
\email{tevfiksah@gmail.com}
\subjclass{53A35, 53A55, 49Q20}
\keywords{Galilean space, relaxed elastic line of second kind, Variational
method, Geodesic}

\begin{abstract}
In this paper, we consider the classical variational problem in the Galilean space. we develop the Euler-Lagrange equations for a elastic line on an oriented surface in the Galilean 3-dimensional space $G_{3}.$ Using the variation method, we will try to give some characterization for the solution curve (the elastic line) of energy equation described by the total squared curvature function of a curve on an oriented surface in $G_3$. Finally, we will investigate whether or not the relaxed elastic curves are on a geodesic. 

\end{abstract}

\maketitle

\section{ Introduction}

There has been in recent years a wide variety of applications of variational methods to various fields of mechanics and technology. Calculus of variation is an important theory which has common applications in geometry, analysis, physics, chemistry and engineering. Variation method is used to find the maxima or the minima of expressions involving unknown functions called functionals. Therefore, the main problem in the calculus of variations is to minimize (or maximize) not only functions but also functionals. For instance; a typical variational problem is to find the fixed length $l$ of a curve joining two points on a surface. The Euler-Lagrange equations is used for the solution of the problem.

In this paper, we will first define the total squared curvature function for a curve in a Galilean space. Using the variation method, we will try to give some characterization for the solution curve (the relaxed elastic line) of energy equation described by the total squared curvature function of a curve on an oriented surface in Galilean space. Finally, we will investigate whether or not the relaxed elastic curves are on a geodesic. 

 There has been in recent years a wide variety of applications of variational methods to various geometries. For instance, several geometers were
interested in studying of calculus of variations in Euclidean Space
\cite{manning,nickerson,unan,yilmaz}. Nevertheless, similar applications of this theory in
Minkowski space can be found in \cite{tutar}. The main point of the studying the
variation calculus
 is defining the Euler-Lagrange equations for a relaxed elastic
line on an oriented surface.

The natural variational integrals in geometry are the common integrals on
space curves $\alpha \left( s\right) $. These include the length $L\left(
\alpha \right) =\dint ds$, total squared curvature $K\left( \alpha \right)
=\dint \kappa ^{2}ds$ used in \cite{manning, nickerson}, total squared torsion $T\left( \alpha
\right) =\dint \tau ^{2}ds$ used in \cite{tutar,unan} and the integral $H\left( \alpha
\right) =\dint \kappa ^{2}\tau ds$ used in \cite{sahin}. An elastic line of length $l$
is defined as a curve with associated energy equation $K=\dint%
\limits_{0}^{l}\kappa ^{2}\left( s\right) ds$ by Nickerson and Manning in
\cite{nickerson}, where $s$ is the arc-length along the curve, $\kappa ^{2}\left(
s\right) $ is the square curvature in there. The integral $K$ is called the
total square curvature. If no boundary conditions are imposed at $s=l$, and
if no external forces act at any $s,$ the elastic line is \textit{relaxed,}
[manning,nickerson]. 
Hilbert and Cohn-Vossen stated in \cite{hilbert} that a relaxed elastic line with
specified position and tangent at $s=0$ always has the trajectory of a
geodesic. However, Manning have proved in \cite{manning} that the
conslusion of Hilbert and Cohn-Vossen is incorrect. Manning in
\cite{manning} has obtained Euler-Lagrange equations of the relaxed elastic line on an
oriented surface as a model of DNA molecule. 

Besides Euclidean Geometry, a range of new types of geometries have been
invented and developed in the last two centuries. They can be introduced in
a variety of ways. One possible way is through projective manner, where one
can express metric properties through projective relations. Among these
geometries, there is also Galilean geometry which is our matter in this
paper.

In literature, there is very ? study on relaxed elastic line on an oriented surface in Galilean space $
G_{3}$. Gokce and Sahin in \cite{gokce,sahin} have obtained intrinsic equations of the relaxed elastic line on an oriented surface in $G_3$ for different energy functions. However,  Gokce and Sahin in \cite{gokce,sahin} have not use variation method. Hence, in
literature, there is no any study that using of variation method on the relaxed elastic line on an oriented surface in $G_{3}$. 

Therefore, firstly, we will derive the Euler-Lagrange
equations for a relaxed elastic line on an oriented surface in $G_{3}$.  
Hence, we have proved in this paper that the conslusion of Hilbert and Cohn-Vossen is incorrect in $G_3$.  
The main result of this paper is Theorem $2$. We describe the geometric
meaning of Theorem $2$ in section 4.

\section{Preliminaries on Galilean Geometry}

\qquad \textquotedblright All geometry is projective
geometry.\textquotedblright\ (A. Cayley). From A. Cayley point of view, $G_{3
\text{ }}$ is a real 3-dimensional projective space $P^{3}(\mathbb{R}) $, is the set of equivalence classes of $\sim $ on $\mathbb{R}^{4}-\left\{ 0\right\} $ by equivalence relation $x\sim y$ iff 
$x\lambda y$
for some $\lambda \in \mathbb{R}\backslash \{ 0\} $. Thus, $P^{3}(\mathbb{R}) $ obtained as a factor space on $\mathbb{R}^{4}\backslash \{ 0\} $ by $\sim $, i.e. $P^{3}(\mathbb{R}) \widetilde{=}( \mathbb{R}^{4}-\{ 0\} ) /\sim $ \cite{cox}. We can think of $P^{3}(\mathbb{R}) $ more geometrically as a set of lines through the origin in 
$\mathbb{R}^{4}$. $G_{3}$ is a real Cayley-Klein space equipped with the projective
metric of signature $( 0,0,+,+) $, as showed in \cite{molnar}. The
absolute of the Galilean geometry is an ordered triple $\{w,f,I\}$, where $w$
is the ideal (absolute) plane, $f$ \ is the line (absolute line) in $w$ and $
I$ is the fixed elliptic involution of points of $f$. \ The points, the
lines and the planes of $P^{3}( \mathbb{R}) $ are the one-dimensional, two-dimensional and three-dimensional
subspaces of $\mathbb{R}^{4}$, respectively \cite{casse}. Therefore, $G_{3}$ contains $\mathbb{R}^{3}$ as a proper subset and the complement in $G_{3}$ to $w$ is
diffeomorphic to $\mathbb{R}^{3}$.

Let $P$ be any point of $\mathbb{R}^{3}$ with affine coordinates $( x,y,z) $. Write $(
x,y,z) $ as $( \frac{X_{1}}{X_{0}},\frac{X_{2}}{X_{0}},\frac{X_{3}
}{X_{0}}) $, where $X_{0}$ is some common deminator. Call $(
X_{0},X_{1},X_{2},X_{3}) $ the \textit{homogeneous coordinates of }$P$
. Thus, the homogeneous coordinates $( X_{0}:X_{1}:X_{2}:X_{3}) $
and $\rho ( X_{0}:X_{1}:X_{2}:X_{3}) $ refer to the same point
for all $\rho \in \mathbb{R}-\left\{ 0\right\} ,$ \cite{casse}. Now,we can introduce homogeneous coordinates in
$G_{3}$ in such a way that the absolute plane $w$ is given by $X_{0}=0$, the
absolute line $f$ by $X_{0}=X_{1}=0$ and the elliptic involution $I$ by
\begin{equation*}
( 0:0:X_{2}:X_{3}) \rightarrow ( 0:0:X_{3}:-X_{2}) .
\end{equation*}

In the nonhomogeneous coordinates the isometries group $B_{6}$ has the form
\begin{eqnarray}
\overline{x} &=&a+x  \notag \\
\overline{y} &=&b+cx+y\cos \varphi +z\sin \varphi \\
\overline{z} &=&d+ex-y\sin \varphi +z\cos \varphi  \notag
\end{eqnarray}
where $a,b,c,d,e$ and $\varphi $ are real numbers. The group of motions of $%
G_{3}$ is a six-parameter group \cite{pavkovic}.\newline
The Galilean norm of the vector $\textbf{v}=(x,y,z)$ defined by
\begin{equation}
\left\Vert \mathbf{v}\right\Vert _{G}=%
\begin{Bmatrix}
x\ \ ,\  & if\ \ \ x\neq 0 \\ 
\sqrt{ y^{2}+z^{2} }, & if\ \ \ x=0\ 
\end{Bmatrix}
.
\end{equation}
A vector $\mathbf{v}=( x,y,z) $ is said to be non-isotropic if $x\neq 0$.
All unit non-isotropic vectors are of the form $\left( 1,y,z\right) $.
For a curve $\alpha :I\rightarrow G_{3}$ , \ $I\subset \mathbb{R}$ parametrized by the arc-length parameter $s=x$, given in the coordinate form
\begin{equation}
\alpha \left( x\right) =\left( x,y\left( x\right) ,z\left( x\right) \right) ,
\end{equation}
the curvature $\kappa \left( x\right) $ and the torsion $\tau \left(
x\right) $ are defined by
\begin{eqnarray}
\kappa \left( x\right) &=&\Vert \alpha''(x) \Vert_{G}=\sqrt{y^{\prime \prime }\left( x\right)
^{2}+z^{\prime \prime }\left( x\right) ^{2}}, \\
\tau \left( x\right) &=&\frac{\det \left( \alpha ^{\prime }\left( x\right)
,\alpha ^{\prime \prime }\left( x\right) ,\alpha ^{\prime \prime \prime
}\left( x\right) \right) }{\kappa ^{2}\left( x\right) }
\end{eqnarray}
and the associated moving trihedron is given by
\begin{eqnarray}
T\left( x\right) &=&\alpha ^{\prime }\left( x\right) ,  \notag \\
N\left( x\right) &=&\frac{\alpha''(x)}{\kappa \left( x\right) } , \\
B\left( x\right) &=&\frac{1}{\kappa \left( x\right) }\left( 0,-z^{\prime
\prime }\left( x\right) ,y^{\prime \prime }\left( x\right) \right) .  \notag
\end{eqnarray}
The vectors $T\left( x\right) ,N\left( x\right) $ and $B\left( x\right) $
are called the vectors of the tangent, principal normal and the binormal
line, respectively \cite{pavkovic,sahinac}. Therefore, the Frenet-Serret formulas can be
written in matrix notation as
\begin{equation}
\begin{bmatrix}
T \\ 
N \\ 
B
\end{bmatrix}
^{\prime }=
\begin{bmatrix}
0 & \kappa & 0 \\ 
0 & 0 & \tau \\ 
0 & -\tau & 0
\end{bmatrix}
\begin{bmatrix}
T \\ 
N \\ 
B
\end{bmatrix}.
\end{equation}
Let $\mathbf{a=}( a_{1},a_{2},a_{3}) $, $\mathbf{b=}(b_{1},b_{2},b_{3}) $ be two vectors in $G_{3}.$

For any regular curve $\alpha $ the following formulas hold:
\begin{equation}
\kappa =\frac{\left\Vert \gamma ^{\prime }\times _{G}\gamma ^{\prime \prime
}\right\Vert_{G} }{\left\Vert \gamma ^{\prime }\right\Vert ^{3}_{G}},\text{ \ }\tau =
\frac{\det \left( \gamma ^{\prime },\gamma ^{\prime \prime },\gamma ^{\prime
\prime \prime }\right) }{\left\Vert \gamma ^{\prime }\times _{G}\gamma
^{\prime \prime }\right\Vert ^{2}_{G}}.
\end{equation}
Here,the Galilean cross product $\times _{G}$ is defined by
\begin{equation}
\mathbf{a}\times _{G}\mathbf{b=}
\begin{vmatrix}
0 & e_{2} & e_{3} \\ 
a_{1} & a_{2} & a_{3} \\ 
b_{1} & b_{2} & b_{3}
\end{vmatrix}
\end{equation}
as in \cite{sipus}.

The Galilean Sphere $S_{G}^{2}$ is defined by $S_{G}^{2}=\{ (
x,y,z) \in G_{3}~\ :~\vert x-x_{0}\vert =r\} .$
\newline
For more facts about the Galilean geometry, we refer the reader to
\cite{gokce, pavkovic, sahinuk, sahin, roschel, yaglom, sipus} and references therein.

\begin{theorem}
Let $S$ be a surface in $G_{3}$ and $\alpha $ be a curve on $S$. At a point 
$\alpha (x)$ of $\alpha$, let $T$ denote
the unit tangent vector of $\alpha$ at $\alpha(x)$, $n$ the unit normal to $S$ and $
n\times _{G}T=Q$ the tangential-normal. Then $\{ T,Q,n\} $ is an
orthonormal basis at $\alpha(x)$ in $S$. This frame is called Galilean Darboux frame or tangent-normal frame. The Frenet-Serret formulas the for the Galilean Darboux frame can be written in matrix notation as
\begin{equation}
\left[ 
\begin{array}{c}
T \\ 
Q \\ 
n
\end{array}
\right] ^{\prime }=\left[ 
\begin{array}{ccc}
0 & \kappa _{g} & \kappa _{n} \\ 
0 & 0 & \tau _{g} \\ 
0 & -\tau _{g} & 0
\end{array}
\right] \left[ 
\begin{array}{c}
T \\ 
Q \\ 
n
\end{array}
\right]
\end{equation}
where $\kappa _{g}$,$\kappa _{n},\tau _{g}$ are geodesic curvature, normal
curvature and geodesic torsion, respectively.
\end{theorem}

Eq.(2.10) implies the important relation $\kappa ^{2}=T^{\prime }._{G}$ $
T^{\prime }=\kappa _{g}^{2}+\kappa _{n}^{2}$ where $\kappa ^{2}$ is the
square curvature of $\alpha $.\newline

\begin{proof} It is clear?
$( 2.7) $ and $(2.10)$ implies the important relation.
\end{proof}

We shall treat a $C^{3}$- surface, $r\geq 1$, as a subset $S\subset G_{3}$
for which there exists an open subset $D$ of $\mathbb{R}^{2}$ and a $C^{3}$- mapping $\psi :D\rightarrow G_{3}$ satisfying $S=\psi\left( D\right) $. A $C^{3}$- surface $S\subset G_{3}$ is called regular if $\psi $ is an immersion, and simple if $\psi $ is an embedding.
Therefore, the oriented surface $S$ represented by $\psi (u_{1},u_{2}) =( \mathbf{X}( u_{1},u_{2}) ,\mathbf{Y}( u_{1},u_{2}) ,\mathbf{Z}( u_{1},u_{2}) ) $.
Here $\mathbf{X,Y}$ and $\mathbf{Z}$ are the coordinate functions of $\psi$
. A unit normal field of surface at a point $p$ is defined by
\begin{equation}
n=\frac{\psi _{1}\times _{G}\psi _{2}}{\left\Vert \psi _{1}\times _{G}\psi
_{2}\right\Vert }=\frac{1}{W}\left( 0,\mathbf{X}_{1}\mathbf{Z}_{2}-\mathbf{X}
_{2}\mathbf{Z}_{1},\mathbf{X}_{2}\mathbf{Y}_{1}-\mathbf{X}_{1}\mathbf{Y}
_{2}\right) ,
\end{equation}
\newline
where $\psi _{i}=\frac{\partial \psi }{\partial u_{i}}, W=\Vert \psi
_{u_{1}}\times _{G}\psi _{u_{2}}\Vert, \, \mathbf{X}_{i}=\frac{\partial \mathbf{X}}{
\partial u_{i}}, \mathbf{Y}_{i}=\frac{\partial \mathbf{Y}}{\partial u_{i}}$ and $\mathbf{Z}_{i}=\frac{
\partial \mathbf{Z}}{\partial u_{i}}$, $i=1,2.$ In a tangent plane of the surface at the
point $p$, there is a unique isotropic direction defined by the condition $
\mathbf{X}_{1}du_{1}+\mathbf{X}_{2}du_{2}=0$. A side tangential vector $Q=
\frac{1}{W}( \mathbf{X}_{2}\psi _{1}-\mathbf{X}_{1}\psi _{2}) $ is a unit isotropic
vector in a tangent plane. The curve $\alpha $ on the oriented surface $S$
is specified by $\psi(x) =\psi( u_{1}(x), u_{2}(x)).$ Then the unit tangent vector along $\alpha $
is
\begin{equation}
T=\frac{d\psi }{dx}=\psi _{1}\overset{\centerdot }{u_{1}}+\psi _{2}\overset{
\centerdot }{u_{2}},
\end{equation}
where $\overset{\centerdot }{u_{1}}=\frac{du_{1}}{dx}$ and $\overset{
\centerdot }{u_{2}}=\frac{du_{2}}{dx}$. \newline
The first fundamental form of a surface is introduced in the following way
\begin{equation}
dx^{2}=\left( \mathbf{X}_{1}du_{1}+\mathbf{X}_{2}du_{2}\right)
^{2}+\varepsilon \left\{ \left( \mathbf{Y}_{1}du_{1}+\mathbf{Y}
_{2}du_{2}\right) ^{2}+\left( \mathbf{Z}_{1}du_{1}+\mathbf{Z}
_{2}du_{2}\right) ^{2}\right\} ,
\end{equation}
\newline
where
\begin{equation}
\varepsilon =\left\{ 
\begin{array}{c}
0,\text{ \ \ }dx\text{\ is non-isotropic,} \\ 
1,\text{ \ }dx\text{ \ \ \ \ \ is isotropic }
\end{array}
\right\} .
\end{equation}
\newline
Since $T.T=1,$ we have a constraining relation between any pair of functions 
$\left( u_{1}\left( x\right) ,u_{2}\left( x\right) \right) $ that define the
curve $\alpha $ on the surface $S:$
\begin{equation}
g\left( u_{1},u_{2},\overset{\centerdot }{u_{1}},\overset{\centerdot }{u_{2}}
\right) =1,
\end{equation}
\newline
where
\begin{equation}
g=\tsum\limits_{i,j=1}^{2}g_{ij}\dot{u}_{i}\dot{u}_{j},
\end{equation}
\newline
and $g_{1}=\mathbf{X}_{1},$ $g_{2}=\mathbf{X}_{2}$ and $g_{ij}=g_{i}.g_{j},$ 
$i,j=1,2,$ where $g_{ij},$ $i,j=1,2,$ are the coefficients of the first
fundamental form. The normal curvature can be expressed in terms of the
coordinates $\left( u_{1}\left( x\right) ,u_{2}\left( x\right) \right) $
along the curve $\alpha $ as
\begin{equation}
\kappa _{n}=\tsum\limits_{i,j=1}^{2}L_{ij}u_{i}^{\prime }u_{j}^{\prime },
\end{equation}
where $L_{11},L_{12}$ and $L_{22}$ are the coefficients of the second
fundamental form $II$ and the normal components of $\psi _{11},\psi _{12}$
and $\psi _{22}$, respectively. It holds
\begin{equation}
L_{ij}=\frac{\mathbf{X}_{2}\psi _{ij}-\mathbf{X}_{ij}\psi _{2}}{\mathbf{X}
_{2}}\cdot n\text{ or }L_{ij}=\frac{\mathbf{X}_{1}\psi _{ij}-\mathbf{X}
_{ij}\psi _{1}}{\mathbf{X}_{1}}\cdot n.
\end{equation}
\newline
Similarly, for the geodesic torsion $\tau _{g}$ we have
\begin{equation}
\tau _{g}=\tsum\limits_{i,j=1}^{2}g^{i}L_{ij}u_{i}^{\prime },
\end{equation}
\newline
where $g^{1}=\frac{\mathbf{X}_{2}}{W}$, $g^{2}=-\frac{\mathbf{X}_{1}}{W}$
and $g^{ij}=g^{i}.g^{j}.$ On the other hand, the square geodesic curvature can
be obtained as
\begin{equation}
\kappa _{g}^{2}=\tsum\limits_{i,j=1}^{2}g_{ij}\gamma _{i}\gamma _{j},
\end{equation}
\newline
where
\begin{equation}
\gamma _{i}=\overset{\cdot \cdot }{u}_{i}+\tsum\limits_{k,l=1}^{2}\Gamma
_{kl}^{i}\dot{u}_{k}\dot{u}_{l},i=1,2
\end{equation}
\newline
and the quantities $\Gamma _{kl}^{i}$ are the Christoffel symbols of the
second kind; available formulas are expressed as functions $g_{ij}$, and
their first partial derivatives with respect to $u_{i}$. Furthermore, the
Christoffel symbols can be written as follows:
\begin{equation}
\Gamma _{kl}^{1}=\frac{\mathbf{X}_{2}\psi _{ij}-\mathbf{X}_{ij}\psi _{2}}{W}
\cdot Q,\text{ }\Gamma _{kl}^{2}=\frac{\mathbf{X}_{1}\psi _{ij}-\mathbf{X}
_{ij}\psi _{1}}{W}\cdot Q,\text{ }k,l=1,2.
\end{equation}
\newline
Hence, the equations of a geodesic curve, which is characterized by
identically vanishing $\kappa _{g}$, must be given by $\gamma _{i}=0,$ and
indeed they are \cite{carmo}.



\section{The incomplete and complete variational problems}

Now we can give the definition of elastic line on an oriented surface $S$ in $G_{3}.$\newline
\subsection{The incomplete variational problem}
\begin{definition}
Let $\alpha $ be a $C^{2}$-curve with parametrized by arc-length $x,$ $0\leq
x\leq \ell ,$ on an oriented surface $S$ in $G_{3}.$ An elastic line of length $\ell $ is defined as a curve with associated energy
\begin{equation}
K=\dint\limits_{0}^{\ell }\kappa ^{2}dx,
\end{equation}
where $\kappa ^{2}$ is the square curvature of the curve $\alpha$ \cite{manning}.
The integral $K$ is called the \textit{total square curvature}. Therefore, the arc $\alpha$ is called an elastic line if it is an extremal for the variational problem of minimizing the value of $K$ with in the family of all arcs of length $\ell $ on $S$, having the same initial point and initial direction as $\alpha .$
\end{definition}

The pertinent problem is to identify that curve in a family of curves of
fixed length $\ell $ which has the least total square curvature. In general the minimizing curve is not a geodesic line. It is a geodesic line, however, in the circumstances described below.\newline
Let
\begin{equation}
K_{n}=\dint\limits_{0}^{\ell } \kappa^2 _{n}(x)dx.
\end{equation}
\newline
and let $\alpha $ be the curve that minimizes $K_{n}$ among all curves of
length $\ell $ on a surface $S$. \newline
For any curve $\alpha $ that minimizes $K_n$ on $S$, if $\alpha $ is geodesic line, then the
geodesic curvature $\kappa _{g}$ identically vanishes. Hence $K=K_{n}$. We call this problem \textit{incomplete}
because it seeks to minimize $K_{n}$, not the total square curvature $K$. If a curve that minimises $K_n$ is a geodesic, then the curve is minimise $K$. Therefore, the relaxed elastic curve lie on a geodesic trajectory.
\newline
Through $\kappa_n^2$ has a dependence
on $u_{1},u_{2},\dot{u}_{1}$ and $\dot{u}_{2}$. Therefore, the Euler-Lagrange equations for the
incomplete variational problem are
\begin{eqnarray}
H_{u_{1}}-( H_{\dot{u}_{1}}) ^{\cdot } &=&0, \\
H_{u_{2}}-( H_{\dot{u}_{2}}) ^{\cdot } &=&0,
\end{eqnarray}
\newline
where
\begin{equation}
H=\kappa _{n}^{2}+\lambda ( g-1)
\end{equation}
\newline
and $\lambda =\lambda \left( x\right) $ is a Lagrange multiplier function.
Equations $(3.3) ,(3.4) $ and $(3.5) $ are
a system of three equations to determine the three functions $u_{1}(x)$, $u_{2}(x)$ and $\lambda(x)$.
Equations $(3.3, 3.4) $ and $(3.5) $ are second order in $u_{1},$ $u_{2}$ and first order in $\lambda$. Therefore, the system is fifth order. To get the system in normal form, Eq. $(2,15)$
must be differentiated one times with respect to $x$. Reintegration of $g'=0$ gives that the constant of integration equal unity. The other four constants of integrations in the general solution are determined by the boundary conditions.\newline
We turn now to the boundary conditions at $x=0$ and $x=\ell$ that determine
the other  four constants of integration of the equations.
This conditions are
\begin{eqnarray}
H_{\dot{u}_1}=H_{\dot{u}_2} =0, \, (x=0, \ell)
\end{eqnarray}
If the initial point of the elastic line is fixed say,
\begin{equation}
u_1(0)=u_{1_0}, \, u_2(0)=u_{2_0}
\end{equation}
but the end of the elastic line is free, then
\begin{equation}
H_{\dot{u}_1}=H_{\dot{u}_2} =0, \, (x=\ell).
\end{equation}
Therefore, the four constants of integration are now determined by above two equations.

\subsection{The complete variational problem} 
If a curve that minimises $K_n$ is not a geodesic, then the curve is or not minimise $K$. Therefore, the incomplete variational problem is not anything say about the relaxed elastic curve that minimises K. In this case, we have no choice but is to find the curve that minimises the integral
$$K=\dint\limits_{0}^{\ell }\kappa ^{2}dx=\dint\limits_{0}^{\ell }(\kappa_n ^{2}+\kappa_g^2)dx$$
The relaxed elastic curve that is not a geodesic curve is one of solutions of the complete variational problem.
Through $\kappa_n^2+\kappa_g^2$ has a dependence
on $u_{1},u_{2},\dot{u}_{1}, \dot{u}_{2}, \ddot{u}_{1}$ and $\ddot{u}_{2}$. Therefore, the Euler-Lagrange equations for the complete variational problem are
\begin{eqnarray}
H_{u_{1}}-( H_{\dot{u}_{1}})^{\cdot}+ ( H_{\ddot{u}_{1}})^{\cdot\cdot}&=&0, \\
H_{u_{2}}-( H_{\dot{u}_{2}}) ^{\cdot }+( H_{\ddot{u}_{2}})^{\cdot\cdot} &=&0,
\end{eqnarray}
Now, let
\begin{equation}
H=\kappa^{2}+\lambda \left( g-1\right)
\end{equation}
with the constrain $g=1$ given by Eqs. $(2.15) $, $(
2.16) $ and $\lambda =\lambda(x)$ again a Lagrange
multiplier function to be determined in the course of solving the system of
differential equations, then we call this variation problem \textit{complete}
because it seeks to minimize the total square curvature.
Equations $(3.9), (3.10) $ and $(3.11)$ are
a system of three equations to determine the three functions $u_{1}(x) $, $u_{2}(x)$ and $\lambda(x)$.
Equations $(3.9)$ and $(3.10)$ are fourth  order in $u_{1},$ $u_{2}$ and first order in $\lambda $. Therefore, the system is of order nine. To get the system in normal form, Eq. $(2.15)$
must be differentiated three times with respect to $x$. Reintegration of $g'''=0$ gives the respective fixed values $0, 0,$ and $1$ to the resulting three constants of integration. The other six constants of integrations in the general solution are determined by the boundary conditions.
Because the side condition Eq. (2.15) does not involve $\ddot u_1$ and $\ddot u_2$, the boundary terms of the variation are not determined in the usual straightforward way. We nevertheless find the following six boundary terms if Eq. (2.15) can be solved for $\dot u_2$ as a function of $u_1, u_2,$ and $\dot u_1$ to give $\dot u_2= U_2(u_1,u_2,\dot u_1)$:
\begin{eqnarray}
\big\{H_{\ddot u_1}+{H_{\ddot u_2}}{(U_2)}_{\dot u_1}\big\} \delta{\dot u_1}    &\hskip 1cm (s=0, \ell)& \\
\big\{H_{\dot{u}_{1}}- ( H_{\ddot{u}_{1}})^{\cdot}+{H_{\ddot u_2}}{(U_2)}_{u_1}\big\} \delta{u_1} &\hskip 1cm (s=0, \ell)& \\
\big\{H_{\dot{u}_{2}}- ( H_{\ddot{u}_{2}})^{\cdot}+{H_{\ddot u_2}}{(U_2)}_{u_2}\big\} \delta{u_2} &\hskip 1cm (s=0, \ell)&
\end{eqnarray}
On the other hand, if ${U_2}_{\dot u_1}$ is singular, then we must use a different set of six boundary terms, after solving Eq. (2.15) for $\dot u_1= U_1(u_1,u_2,\dot u_2)$:
\begin{eqnarray}
\big\{H_{\ddot u_2}+{H_{\ddot u_1}}{(U_1)}_{\dot u_2}\big\}\delta{\dot u_2}    &\hskip 1cm (s=0, \ell)& \\
\big\{H_{\dot{u}_{2}}- ( H_{\ddot{u}_{2}})^{\cdot}+{H_{\ddot u_1}}{(U_1)}_{u_2}\big\} \delta{u_2} &\hskip 1cm (s=0, \ell)& \\
\big\{H_{\dot{u}_{1}}- ( H_{\ddot{u}_{1}})^{\cdot}+{H_{\ddot u_1}}{(U_1)}_{u_1}\big\} \delta{u_1} &\hskip 1cm (s=0, \ell)&
\end{eqnarray}
Suppose that we place no restrictions at all on the elastic line other than confinement to the surface. Then six integration constants are determined by setting equal to zero each of the six factors multiplying $\delta u_1, \delta u_2,$ and $\delta \dot{u_1}$ in the "$U_1$ boundary terms," or those multiplying $\delta u_1, \delta u_2,$ and $\delta \dot{u_2}$ in the "$U_2$ boundary terms." These boundary conditions are \textit{completely natural (free)} \cite{manning}.

If natural conditions are allowed to prevail at $x=\ell$, but the initial position and direction of the elastic line are specified, then we get a different set of equations to determine the six integration constants. We illustrate with the $U_2$ boundary terms. Three boundary conditions are obtained by setting equal to zero each of the three factors multiplying $\delta u_1, \delta u_2,$ and $\delta \dot{u_1}$, with the factors evaluated at $x=\ell$. Two more boundary conditions are specified by the values of $u_1$ and $u_2$ at $x=0$: ${u_1}(0)={u_1}_0, \, {u_2}(0)={u_2}_0$. The unit tangent vector $T(x)$ at any point $x$ of the curve is given by $\phi_{u_1}(u_1,u_2){\dot {u_1}(x)} +\phi_{u_2}(u_1,u_2){\dot {u_2}(x)}$, and so the initial tangent $T(0)=\phi_{u_1}({u_1}_0,{u_2}_0){\dot {u_1}(0)} +\phi_{u_2}({u_1}_0,{u_2}_0){\dot {u_2}(0)}$. Specifications of $\dot{u_1}(0)=\dot{u_1}_0$ therefore completely determines the initial direction, since $\dot{u_2}(0)$ is then fixed at the value $U_2({u_1}_0,{u_2}_0, \dot{u_1}_0)$. Thus, the required sixth boundary condition is that $\dot{u_1}(0)$ have a specified value.

\section{Elastic Lines On Some Surfaces in Galilean Space}

\subsection{Elastic lines on Galilean plane}

Let $S$ be a Galilean plane. Then, for the coefficients of the first and
second fundamental forms we have
\begin{eqnarray*}
g_{11} &=&1,g_{12}=0,g_{22}=1 \\
L_{11} &=&L_{12}=L_{22}=0.
\end{eqnarray*}
Thus, we obtain $\Gamma _{kl}^{j}=0,1\leq i,k,l\leq 2.$ Therefore, we can
say that the normal curvature $\kappa _{n}$ and the geodesic torsion $\tau
_{g}$ of any curve on $S$ vanish in all points of this curve from Eq.(2.17). In particular, any geodesic line on the plane minimizes $K_n$ (giving to it the value zero), Hence, the curve minimizes K. Therefore, the relaxed elastic lines on Galilean plane lies along a geodesic line.

It follows from Eq.(2.21) that the differential equations of geodesic curves on a plane become $u_1''=u_2''=0$. The solutions of differential equations are straight lines $u_1=ax+b, u_2=cx+d, a^2+c^2=1$. 
The relaxed elastic line on a plane assumes the form of a straight line.

\subsection{Elastic lines on Galilean sphere}

Non-isotropic Galilean unit sphere is defined a set of non-isotropic unit
vectors in Galilean space. Therefore, Non-isotropic Galilean unit sphere is
described as set of points $\left( \pm 1,y,z\right) ,$ i.e. is two parallel
plane given by $x=\pm 1.$ The non-isotropic Galilean unit sphere with
parametrization is obtained by $\varphi (u_1,u_2)=( \pm 1,u_1,u_2).$ Its
first fundamental form is given by 
\begin{equation*}
ds^{2}=du_1^{2}+du_2^{2}.
\end{equation*}
Then, the coefficients of the first and second fundamental forms are
\begin{eqnarray*}
g_{11} &=&1,g_{12}=0,g_{22}=1, \\
L_{11} &=&L_{12}=L_{22}=0.
\end{eqnarray*}
Thus, we obtain $\Gamma _{kl}^{j}=0,1\leq i,k,l\leq 2.$ 

Therefore, we can say that the relaxed elastic line on the Galilean sphere assumes the form of a straight line.

\subsection{Elastic lines on a cylinder}
Let $M$ be the cylinder over a curve $C: f(y,z)=c$ in the $yz$ plane. If $\alpha=(0, \alpha_1,\alpha_2)$ is a parametrization of $C$, we assert that $$\chi(u,v)=(u, \alpha_1(v),\alpha_2(v))$$
is a parametrization of $M$. If $\alpha=(0, R\cos(v/R),R\sin(v/R))$ is a parametrization of $C$, we have  $$\chi(u,v)=(u, R\cos(v/R),R\sin(v/R))$$

From Eq.(2.11-2.22), we obtain the following coefficients
\begin{eqnarray*}
g_{11} &=&1,g_{12}=g_{22}=0, \\
L_{11} &=&L_{12}=0, L_{22}=R^{-1}, \\
\Gamma _{kl}^{j}&=&0.
\end{eqnarray*}
where $1\leq i,k,l\leq 2.$  From Eq.(2.20), (2.21) and (2.22) the geodesic curvature along any curve $[u(x),v(x)]$ on cylinder is given by $\kappa_g=\vert  \overset{..}{u}\vert$. It is a characteristic of geodesic curves that the geodesic curvature $\kappa_g$ vanishes identically along them. Therefore, the geodesic on a cylinder have $u(x)=ax+b, \, a^2=1$. Consequently, any curve $[u(x),v(x)]$ on cylinder is geodesic if and only if $u(x)=\pm x+b$.

From Eq.(2.17) the square normal curvature along any curve $[u(x),v(x)]$ on the cylinder is given by 
\begin{equation}
\kappa_N^2(x)=R^{-1}\overset{.}{v}^2
\end{equation}
We obtain following $H$  for the incomplete variational problem from Eq.(2.16) and (3.5)
$$H=R^{-2}\overset{.}{v}^4+\lambda(\overset{.}{u}^2-1).$$
Therefore, we get from the Euler equations (3.3) and (3.4), and also (3.6)
$$H_{\dot{u}}= H_{\dot{v}}=0.$$
Therefore, we have the system of three differential equations for $u(x), v(x),$ and $\lambda(x)$ is 
\begin{eqnarray}
2\lambda{\dot{u}}&=&0, \\
4R^{-2}{\dot{v}^3}&=&0, \\
\dot{u}^2&=&1,
\end{eqnarray}
where Eq.(4.4) is side condition Eq.(2.15) for any curve on a cylinder.
Hence, any solution of this system automatically satisfies the boundary conditions. 
From Eq.(4.4), we have $u(x)=\pm x+b$. We must choose $\lambda=0$ and $v=const.$ to satisfy  Eq.(4.2) and (4.3). Hence, all solutions of this system are any vertical generator $u(x)=\pm x+b$, $v=const.$ together with $\lambda=0$. From Eq.(4.1) the total normal curvature is zero for a vertical generator, and this solution provides a minimum for $K_N.$ In the same time, the vertical generator is a geodesic on cylinder. Consequently, we can say that minimum-energy trajectory on a cylinder adopts the form of a vertical generator. That is,  a relaxed elastic line on cylinder is lying on a geodesic.

Now, we examine the complete variational problem. In this case, we have following $H$  for the complete variational problem from Eq.(3.11)
$$H=R^{-2}\overset{.}{v}^4+ \overset{..}{u}^2+\lambda(\overset{.}{u}^2-1).$$
Therefore, we get from the Euler equations (3.9) and (3.10)
\begin{eqnarray*}
H_{\dot{u}}- (H_{\ddot{u}})'&=&const, \\
H_{\dot{v}}- (H_{\ddot{v}})'&=&const.
\end{eqnarray*}

Therefore, we have the system of three differential equations for $u(x), v(x),$ and $\lambda(x)$ is 
\begin{eqnarray}
\lambda{\dot{u}}-\dddot{u}&=&0, \\
4R^{-2}{\dot{v}^3}&=&0, \\
\dot{u}^2&=&1,
\end{eqnarray}

\subsection{Elastic lines on helical surfaces}

Let $\alpha $ be a plane curve and $d$ a given line. A helical surface (or
generalized helicoid) $H$ is a surface obtained when $\alpha $ is displaced
in a rigid screw motion about $d.$ The normal form of the 1-parameter group
of motions called the screw motions in the Galilean space is given by
\begin{eqnarray*}
\overline{x} &=&pt+x, \\
\overline{y} &=&y\cos t+z\sin t, \\
\overline{z} &=&-y\sin t+z\cos t,
\end{eqnarray*}
where $p\in IR.$ By substituting $p=0$, the Euclidean rotations about the $x-
$axis are described, and obtained surface is a surface of revolution.\newline
Since a plane curve $\alpha $ can lie in a Euclidean or in a isotropic
plane, we treat these two cases separately.\newline
Let a plane curve $\alpha $, given by $\gamma (v) =(0,f(v) ,g(v)) ,$ where $f,g\in C^{2}$, be an
admissible curve in a Euclidean plane. If $\alpha $ undergoes the screw
motion, the helical surface $H_{p}$ with parametrization is obtained by
\begin{equation*}
\varphi (u,v) =(pu,f(v) \cos u+g(v) \sin u,-f(v) \sin u+g(v) \cos u) .
\end{equation*}
Its first fundamental form is given by
\begin{equation*}
dx^{2}=( pdu) ^{2}+\epsilon [( f^{2}+g^{2})
du^{2}+2( gf^{\prime }-fg^{\prime }) dudv+( f^{\prime
2}+g^{\prime 2}) dv^{2}] .
\end{equation*}
It is convenient to assume that the curve $\alpha $ is parametrized by the
Euclidean arc length, that is, $f^{\prime 2}+g^{\prime 2}=1.$ Then, the
coefficients of the first and second fundamental forms are%
\begin{eqnarray*}
g_{11} &=&p^{2}, \, g_{12}=0, \, g_{22}=0, \\
L_{11} &=&fg^{\prime }-f^{\prime }g, \, L_{12}=-sgn(p), \, L_{22}=\frac{
g^{\prime \prime }}{f^{\prime }}, \\
\Gamma _{11}^{2} &=&p^{2}\left( -ff^{\prime }-gg^{\prime }\right) 
\end{eqnarray*}
and other the Christoffel symbols are zero. Therefore, we have 
\begin{eqnarray*}
\kappa _{n} &=&(fg^{\prime }-f^{\prime }g)( \overset{.}{u}) ^{2}-2sgn(p) \overset{.}{u}\overset{.}{v}+\frac{g^{\prime \prime }}{f^{\prime }}( \overset{.}{v}) ^{2}, \\
\tau _{g} &=&p\overset{.}{u}-sgn(p) p\frac{g^{\prime \prime }}
{f^{\prime }}\overset{.}{v}, \\
\kappa _{g} &=&\vert p\vert \overset{..}{u}.
\end{eqnarray*}
\newline
Now, let a plane curve $\alpha $ be an admissible curve in an isotropic
plane. We can treat two cases: the first one when the curve $\alpha $ is
parametrized by $\gamma \left( v\right) =\left( g\left( v\right) ,f\left(
v\right) ,0\right) $ and the second when the curve $\alpha $ is parametrized
by $\gamma \left( v\right) =\left( g\left( v\right) ,0,f\left( v\right)
\right) .$ Then the helical surface $H_{i}$ is given by%
\begin{equation*}
\varphi \left( u,v\right) =\left( pu+g\left( v\right) ,f\left( v\right) \sin
u,f\left( v\right) \cos u\right) 
\end{equation*}
and
\begin{equation*}
\varphi \left( u,v\right) =\left( pu+g\left( v\right) ,f\left( v\right) \cos
u,f\left( v\right) \sin u\right) 
\end{equation*}
respectively.\newline
The first fundamental form of a helical surface $H_{i}$ obtained in this way
is given by
\begin{equation*}
dx^{2}=\left( pdu+g^{\prime }\left( v\right) dv\right) ^{2}+\epsilon \left(
f^{2}\left( v\right) du^{2}+f^{\prime 2}\left( v\right) dv^{2}\right) .
\end{equation*}
If the curve $\alpha $ is parametrized by the isotropic arc length, that is
if $g\left( v\right) =v$, then the above form can be simplified. Therefore,
the coefficients of the first and second fundamental forms are
\begin{eqnarray*}
g_{11} &=&p^{2},g_{12}=p,g_{22}=1, \\
L_{11} &=&\frac{f^{2}}{\sqrt{f^{2}+p^{2}f^{\prime 2}}},L_{12}=\frac{
-pf^{\prime 2}}{\sqrt{f^{2}+p^{2}f^{\prime 2}}},L_{22}=\frac{-ff^{\prime
\prime }}{\sqrt{f^{2}+p^{2}f^{\prime 2}}}, \\
\Gamma _{11}^{1} &=&\frac{pff^{\prime }}{f^{2}+p^{2}f^{\prime 2}},\Gamma
_{12}^{1}=\frac{ff^{\prime }}{f^{2}+p^{2}f^{\prime 2}},\Gamma _{22}^{1}=
\frac{-pf^{\prime }f^{\prime \prime }}{f^{2}+p^{2}f^{\prime 2}}, \\
\Gamma _{11}^{2} &=&\frac{-p^{2}ff^{\prime }}{f^{2}+p^{2}f^{\prime 2}}
,\Gamma _{12}^{2}=\frac{-pff^{\prime }}{f^{2}+p^{2}f^{\prime 2}},\Gamma
_{22}^{2}=\frac{p^{2}f^{\prime }f^{\prime \prime }}{f^{2}+p^{2}f^{\prime 2}}.
\end{eqnarray*}
Therefore, we have 
\begin{eqnarray*}
\kappa _{n} &=&\left( \frac{f^{2}}{\sqrt{f^{2}+p^{2}f^{\prime 2}}}\right)
\left( \overset{.}{u}\right) ^{2}+2\frac{pf^{\prime 2}}{\sqrt{
f^{2}+p^{2}f^{\prime 2}}}\overset{.}{u}\overset{.}{v}+\left( \frac{
-ff^{\prime \prime }}{\sqrt{f^{2}+p^{2}f^{\prime 2}}}\right) \left( \overset{
.}{v}\right) ^{2}, \\
\tau _{g} &=&\overset{.}{\left( \frac{f^{2}+p^{2}f^{\prime 2}}{
f^{2}+p^{2}f^{\prime 2}}\right) \overset{.}{u}+}\left( \frac{pff^{\prime
\prime }-pf^{\prime 2}}{f^{2}+p^{2}f^{\prime 2}}\right) \overset{.}{v}, \\
\kappa _{g} &=&\left\vert p\overset{..}{u}+\overset{..}{v}\right\vert .
\end{eqnarray*}


\begin{thebibliography}{99}
\bibitem{casse} R. Casse, Projective Geometry: An Introduction, Oxford Univ.
Press (2006), pp. 45-51.

\bibitem{capo} R. Capovilla, C. Chryssomalakos and J. Guven, Hamiltonians for
curves, J. Phys. A: Math. Gen. 35(2002), 6571-6587.

\bibitem{cox} D. Cox, J. Little and D. O'shea, Ideals, Variets, and Algorithms
(second edition), Springer-Verlag Newyork, (1997), 349-362.

\bibitem{carmo} M. P. Do Carmo, Differential Geometry of Curves and Surfaces,
Prentice-Hall, New Jerse, (1976).

\bibitem{gokce} G. Gökçe, Some Elastic Curves in Galilean Geometry, PhD thesis, Ondokuz Mayıs University, Samsun, 2011.

\bibitem{hilbert} D. Hilbert and S. Cohn-Vossen, Geometry and Imagination, Chelsea,
New York, 1952.

\bibitem{landau} L. D. Landan and E. M. Lifshitz, Theory of Elasticity. Pergamon
Press, Oxford, 1979, 84.

\bibitem{lang} J. Languer, Recursion in Curve Geometry, New York J. Math.
5(1999), 25-51.

\bibitem{manning} G. S. Manning, Relaxed elastic line on a curved surface, Quart.
Appl. Math. 45(3)(1987), 515-527.

\bibitem{molnar} E. Molnar, The projective interpretation of the eight
3-dimensional homogeneous geometries, Beitr. Algebra Geom. 38(1997), 261-288.

\bibitem{nickerson} H. K. Nickerson and G. S. Manning, Intrinsic equations for a
relaxed elastic line on an oriented surface, Geometriae dedicate 27(1988),
127-136.

\bibitem{pavkovic} B. J. Pavkovi\v{c} and I. Kamenarovi\v{c}, The equiform
differential geometry of curves in the Galilean space $G_{3}$, Glas. Mat.
22(42)(1987), 449-457.

\bibitem{sahinuk} T. \c{S}ahin and M. Y\i lmaz,\ On singularities of the Galilean
spherical darboux ruled surface of a space curve in $G_{3}$, Ukrainian Math.
J. 62(10)(2010), 1377-1387.

\bibitem{sahinac} T. \c{S}ahin and M. Y\i lmaz, The rectifying developable and the tangent indicatrix of a curve in Galilean 3-space, Acta Math. Hung., 132(1-2) (2011), 154-167.

\bibitem{sahin} T. \c{S}ahin, Intrinsic Equations for a Generalized Relaxed Elastic Line on an Oriented Surface in the Galilean Space, Acta Mathematica Scientia 33B(3)(2013), 701-711.


\bibitem{unan} Z. \"{U}nan and M. Y\i lmaz, Elastic lines of second kind on an
oriented surface, Ondokuz May\i s \"{U}niv. Fen dergisi\ 8(1)(1997), 1-10.

\bibitem{roschel} O. R\"{o}schel , Die Geometrie des Galileischen Raumes ,
Habilitationssch. , Inst. f\"{u}r Mat . und Angew. Geometrie. (1984).


\bibitem{yaglom} I.M. Yaglom,\ A Simple Non-Euclidean Geometry and Physical Basis,
Spr\i nger-Verlag Newyork, (1979).

\bibitem{yilmaz} M. Y\i lmaz, Some relaxed elastic line on a curved hypersurface,
Pure Appl. Math. Sci. 39(1994), 59-67.

\bibitem{sipus} Zeljka Milin Sipus, Ruled Weingarten surfaces in the Galilean space, Periodica Mathematica Hungarica, 56(2) (2008), 213-225.

\end{thebibliography}
\end{document}